\documentclass[a4paper,11pt]{amsart}

\usepackage{amsthm,amsfonts,amsmath}
\usepackage{amssymb,graphicx}
\usepackage{xcolor}
\usepackage{float}
\usepackage[labelfont={bf},textfont={it}]{caption}
\usepackage[all]{xy}
\usepackage{verbatim}
\usepackage{enumitem}
\usepackage{mathrsfs}
\usepackage{hyperref}

\usepackage{tikz}
\usetikzlibrary{arrows.meta}

\newtheorem{theorem}{Theorem}[section]
\newtheorem*{theorem*}{Theorem}
\newtheorem{lemma}[theorem]{Lemma}
\newtheorem{corollary}[theorem]{Corollary}
\newtheorem{proposition}[theorem]{Proposition}
\theoremstyle{definition}
\newtheorem{definition}[theorem]{Definition}
\newtheorem{example}[theorem]{Example}
\newtheorem{remark}[theorem]{Remark}
\newtheorem{setting}[theorem]{Setting}

\newcommand{\R}{\mathbb{R}}
\newcommand{\C}{\mathbb{C}}
\newcommand{\Sn}{\mathbb{S}}
\newcommand{\HH}{\mathcal{H}}
\newcommand{\II}{\mathrm{II}}
\newcommand{\dd}{\,\mathrm{d}}
\newcommand{\kv}{\vec{\kappa}}
\DeclareMathOperator{\Reg}{Reg}
\DeclareMathOperator{\Sing}{Sing}
\DeclareMathOperator{\Cone}{Cone}

\begin{document}

\title[Gauss--Bonnet theorem for o-minimal surfaces]
{Gauss--Bonnet theorem for o-minimal surfaces with isolated singularities}

\author[L. H. Pinheiro]{Lucas Holanda Pinheiro}
\author[J. E. Sampaio]{Jos\'e Edson Sampaio}
\author[E. C. da Silva]{Euripedes Carvalho da Silva}

\address{Lucas Holanda Pinheiro: Departamento de Matem\'atica, Universidade Federal do Cear\'a, Av. Humberto Monte, s/n Campus do Pici - Bloco 914, 60455-760, Fortaleza-CE, Brazil. E-mail: {\tt lhpholanda@alu.ufc.br}}

\address{Jos\'e Edson Sampaio: Departamento de Matem\'atica, Universidade Federal do Cear\'a, Av. Humberto Monte, s/n Campus do Pici - Bloco 914, 60455-760, Fortaleza-CE, Brazil. E-mail: {\tt edsonsampaio@mat.ufc.br}}

\address{Euripedes Carvalho da Silva: Departamento de Matem\'atica, Instituto Federal de Educa\c{c}\~ao, Ci\^encia e Tecnologia do Cear\'a,
 	      Av. Parque Central, 1315, Distrito Industrial I, 61939-140,
 	      Maracana\'u-CE, Brazil. \newline
 	      and Departamento de Matem\'atica, Universidade Federal do Cear\'a, Rua Campus do Pici, s/n, Bloco 914, Pici, 60440-900,
	      Fortaleza-CE, Brazil.   E-mail: {\tt euripedes.carvalho@ifce.edu.br}}

\thanks{The second named author was partially supported by CNPq-Brazil grants 303375/2025-6, 404611/2025-7 and 405618/2023-9. The third named author was partially supported by CNPq-Brazil grants 400043/2025-4 and 405618/2023-9. This work was supported by the Serrapilheira Institute (grant number Serra -- R-2110-39576).
}


\begin{abstract}
We study the relationship between topology and curvature of surfaces definable in an o-minimal structure. Our main result is a Gauss--Bonnet-type formula for closed (not necessarily compact) definable surfaces with finitely many singular points, expressing the Euler characteristic in terms of the total Gaussian curvature, the local densities at the singular points, and the density at infinity. As consequences, we obtain an inequality for definable minimal surfaces with isolated singularities and, for complex algebraic curves, a formula and an inequality relating the Euler characteristic, the total curvature, the multiplicities at the singular points and the degree.
\end{abstract}

\maketitle

\section{Introduction}

The Gauss--Bonnet theorem is one of the fundamental results of differential geometry: it connects the topology of a compact Riemannian surface with its total Gaussian curvature. For a compact surface $X$ without boundary it states that
\[
\int_X K\dd A=2\pi\chi(X),
\]
where $K$ denotes the Gaussian curvature and $\chi(X)$ is the Euler characteristic of $X$.

For compact singular surfaces, Grieser \cite{Grieser:2003} obtained the following result.

\begin{theorem}[Grieser {\cite[Theorem 1.4]{Grieser:2003}}]\label{thm:grieser}
Let $V\subset\R^N$ be a compact real analytic surface with isolated singular points $p_1,\ldots,p_R$, and endow $V_0=V\setminus\{p_1,\ldots,p_R\}$ with the Riemannian metric induced by a real analytic Riemannian metric on $\R^N$. For $\varepsilon>0$ small, each connected component of a punctured $\varepsilon$-neighbourhood of a singular point is either cone-like or horn-like \cite[Theorem 1.1]{Grieser:2003}; let $C_1,\ldots,C_S$ be the cone-like components (over all singular points) and let $l_i$ be the coefficient of $r$ in the asymptotic expansion \cite[Theorem 1.3]{Grieser:2003} of the length of the curve of points of $C_i$ at distance $r$ from the corresponding singular point. Then $K$ is integrable on $V_0$ and
\[
\chi(V)
=
R+\frac{1}{2\pi}\int_{V_0}K\dd A
-\frac{1}{2\pi}\sum_{i=1}^{S}l_i.
\]
\end{theorem}

Thus horn-like components give no contribution, while each singular point contributes $1-\frac{1}{2\pi}\sum l_i$, the sum being over the cone-like components at that point. When the metric is the Euclidean one, horn-like components have links of length $o(r)$, and $\frac{1}{2\pi}\sum_{i=1}^S l_i=\sum_{j=1}^R\theta(V,p_j)$, where $\theta(V,p_j)$ is the local density (see Subsection~\ref{subsec:density} and Proposition~\ref{prop:local}(b)). An essential ingredient in Grieser's approach is the resolution of singularities, which provides a detailed description of the local geometry and the corresponding asymptotic expansions.

For noncompact surfaces, the relationship between topology and curvature is more subtle, since the geometry at infinity contributes to the global formula. For complete Riemannian surfaces of finite topological type whose total curvature exists, Cohn-Vossen \cite{CohnVossen:1935} proved the inequality $\int_M K\dd A\leq 2\pi\chi(M)$, and Shiohama \cite{Shiohama:1985} identified the defect $2\pi\chi(M)-\int_MK\dd A$ with $\lim_{t\to\infty}L(t)/t=\lim_{t\to\infty}2A(t)/t^2$, where $L(t)$ and $A(t)$ are the length and the area of the intrinsic metric circle and disc of radius $t$. For complete submanifolds of Euclidean space with cone-like ends see Dillen and K\"uhnel \cite{DillenKuhnel:2005}, and for closed semialgebraic sets Dutertre \cite{Dutertre:2008,Dutertre:2012} obtained Gauss--Bonnet formulas involving Lipschitz--Killing curvatures at infinity. In particular, for a smooth closed semialgebraic surface $X\subset\R^n$, \cite[Theorem 4.3]{Dutertre:2012} (case $d=2$) seems to give $\chi(X)=\frac{1}{2\pi}\int_XK\dd A+\theta(X,\infty)$, after identifying the normalizations of the Lipschitz--Killing--Weyl curvatures ($K_0\equiv\operatorname{vol}(\Sn^{n-3})$ and $\frac{1}{s_{n-1}}K_2=\frac{1}{2\pi}K$); see also \cite[Theorem 3.9]{Dutertre:2012} for closed semialgebraic sets with singularities.

O-minimal geometry provides a natural setting for these questions. Definable sets enjoy strong finiteness properties, which make them particularly suitable for the study of singularities and of the asymptotic geometry at infinity. Throughout the paper, a set is said to be \emph{definable} if it is definable in a fixed o-minimal structure on $\R$ (see Definition~\ref{definivel}); we refer to \cite{Coste:2000,Dries:1998} for background. In this setting, the local density and the density at infinity provide quantitative that describe the geometry of a definable set near a point and at infinity, respectively.

In this paper we develop an elementary geometric approach to a Gauss--Bonnet formula for closed definable surfaces with isolated singularities. Our method relies on the analysis of the total geodesic curvature of the links and on its relation with the local density and with the density at infinity (see Subsection~\ref{subsec:density} for the definitions). In contrast with Grieser's approach, we use neither resolution of singularities nor the theory of normal cycles and Lipschitz--Killing curvature measures; the o-minimal finiteness properties are used only through uniform bounds on the number of connected components of definable families. Moreover, our framework applies to surfaces definable in arbitrary o-minimal structures and includes the noncompact case, where the contribution of infinity is expressed in terms of the density at infinity.

For a closed set $X\subset\R^n$ we denote by $\Sing(X)=\Sing_2(X)$ its $C^2$ singular set and by $\Reg(X)=X\setminus\Sing(X)$ its $C^2$ regular part (Definition~\ref{def:sing}); $X$ is \emph{pure two-dimensional} if it has dimension two at each of its regular points (Definition~\ref{def:sing}). We write $\int_XK\dd A:=\int_{\Reg(X)}K\dd A$, and $\chi(X)$ denotes the Euler characteristic of the singular homology of $X$. Our main result is the following.

\begin{theorem}\label{thm:main}
Let $X\subset\R^n$ be a closed, pure two-dimensional set definable in an o-minimal structure, and assume that $\Sing(X)=\{p_1,\ldots,p_s\}$ is finite ($s\geq 0$). Then the homology of $X$ is finitely generated, the Gaussian curvature $K$ of $\Reg(X)$ is integrable, the local densities $\theta(X,p_i)$ and the density at infinity $\theta(X,\infty)$ exist, and
\[
\chi(X)
=
s+\frac{1}{2\pi}\int_X K\dd A
-\sum_{i=1}^{s}\theta(X,p_i)
+\theta(X,\infty).
\]
\end{theorem}

Equivalently, $\chi(X)=\frac{1}{2\pi}\int_XK\dd A+\sum_{i=1}^s\big(1-\theta(X,p_i)\big)+\theta(X,\infty)$: each singular point contributes $1-\theta(X,p_i)$ and infinity contributes $\theta(X,\infty)$. If $X$ is bounded (hence compact), then $\theta(X,\infty)=0$. In the compact real analytic case with the Euclidean metric, Theorem~\ref{thm:main} recovers Theorem~\ref{thm:grieser}.

\begin{example}\label{ex:examples}

(1) For the plane
$X=\{(x,y,z)\in\R^3:z=0\}\cong\R^2$,
we have $\chi(X)=1$, $K\equiv 0$, $s=0$ and
$\theta(X,\infty)=1$.

(2) For the circular cylinder
$X=\{(x,y,z)\in\R^3:x^2+y^2=1\}$,
we have $\chi(X)=0$, $s=0$, $K\equiv 0$ and
$\theta(X,\infty)=0$.

(3) For the upper cone
$X=\{(x,y,z)\in\R^3:z^2=x^2+y^2,\ z\geq0\}$,
we have $\chi(X)=1$, $s=1$ and $K\equiv 0$.
Moreover, $\theta(X,0)=\theta(X,\infty)=1/\sqrt{2}$.

(4) For the double cone
$X=\{(x,y,z)\in\R^3:z^2=x^2+y^2\}$,
we have $s=1$, $\chi(X)=1$ and $K\equiv 0$.
Its link consists of two circles of radius $r/\sqrt{2}$,
hence $\theta(X,0)=\theta(X,\infty)=\sqrt{2}$.

(5) For the paraboloid
$X=\{(x,y,z)\in\R^3:z=x^2+y^2\}$,
we have $\chi(X)=1$, $s=0$ and
$\displaystyle\int_X K\,\dd A=2\pi$,
while $\theta(X,\infty)=0$.

(6) For the catenoid
$X=\{(x,y,z)\in\R^3:x^2+y^2=\cosh^2 z\}$,
which is definable in $\R_{\exp}$, we have
$\chi(X)=0$, $s=0$,
$\displaystyle\int_X K\,\dd A=-4\pi$
and $\theta(X,\infty)=2$ (each of its two ends is asymptotic to a plane).

(7) For the cuspidal surface
$X=\{(x,y,z)\in\R^3:x^2+y^2=z^6,\ z\geq0\}$,
we have $\chi(X)=1$, $s=1$,
$\displaystyle\int_X K\,\dd A=-2\pi$,
$\theta(X,0)=0$ and $\theta(X,\infty)=1$.

(8) If $X$ is a compact surface of revolution which is smooth except at a cuspidal point $p$ where it is locally of the form $\{x^2+y^2=z^{2a},\,z\geq0\}$ with $a>1$, then $\theta(X,p)=0$ and the singular point contributes $1$, as a smooth point would.

\end{example}

\begin{remark}\label{rem:hypotheses}
(1) The hypotheses \emph{closed} and \emph{pure two-dimensional} cannot be removed: the open unit disc $X\subset\R^2$ satisfies $\chi(X)=1$, while the right-hand side of the formula vanishes; and for $X=\R^2\times\{0\}\cup\{(0,0,1)\}\subset\R^3$ (which is closed, with $\Sing(X)=\emptyset$ but with an isolated point) we have $\chi(X)=2$, while the right-hand side equals $1$.

(2) No orientability assumption is needed.

(3) A point $p\in\Sing(X)$ at which $X$ is a $C^1$ submanifold satisfies $\theta(X,p)=1$ and gives no contribution, as expected.

(4) The singular Euler characteristic $\chi(X)$ coincides with the o-minimal Euler characteristic $E(X)$ of \cite[Ch.~4]{Dries:1998}. Indeed, for $R$ large, $X\cap\overline{B_R}$ is a deformation retract of $X$ (Lemma~\ref{lem:structure}), $E(X\cap\overline{B_R})=\chi(X\cap\overline{B_R})$ by definable triangulation \cite[Ch.~8]{Dries:1998}, and $E(X\setminus\overline{B_R})=0$ by Hardt's definable triviality theorem \cite[Ch.~9]{Dries:1998} applied to $x\mapsto|x|$, since the fibres are finite unions of circles.
\end{remark}

A surface $M\subset\R^n$ of class $C^2$ is \emph{minimal} if its mean curvature vector vanishes identically. As a consequence of Theorem~\ref{thm:main} we obtain the following.

\begin{corollary}\label{cor:minimal}
Let $X$ be as in Theorem~\ref{thm:main} and assume that $\Reg(X)$ is a minimal surface. Then
\[
\chi(X)\leq s-\sum_{i=1}^s\theta(X,p_i)+\theta(X,\infty),
\]
and equality holds if and only if $X$ is a finite union of affine $2$-planes.
\end{corollary}

An interesting consequence of Theorem~\ref{thm:main} concerns complex algebraic curves. By a classical result of Thie \cite{Thie:1967} (see also Kurdyka and Raby \cite{KurdykaRaby:1989} for the real subanalytic setting), the local density of a complex analytic curve at a point coincides with its multiplicity. Moreover, the density at infinity of a complex algebraic set coincides with its degree; this goes back to \cite{SilvaSampaio:2026}. Consequently, Theorem~\ref{thm:main} yields the following formula.

\begin{corollary}\label{cor:complex}
Let $X\subset\C^n$ be a complex algebraic curve (i.e., a pure one-dimensional complex algebraic set) such that $\Sing(X)=\{p_1,\ldots,p_s\}$. Then
\[
\chi(X)
=
s+\frac{1}{2\pi}\int_X K\dd A
-\sum_{i=1}^s m(X,p_i)
+\deg(X).
\]
Moreover,
\[
\chi(X)\leq \deg(X)-\sum_{i=1}^s\big(m(X,p_i)-1\big),
\]
and equality holds if and only if $X$ is a union of $\deg(X)$ distinct complex affine lines.
\end{corollary}

For instance, for a smooth affine curve $X$ the corollary gives $\int_XK\dd A=2\pi(\chi(X)-\deg X)$; for the hyperbola $\{xy=1\}\subset\C^2$ this gives $\int_X K\dd A=-4\pi$.

\begin{remark}\label{rem:jorge-meeks}
If $\operatorname{Sing}(X)=\emptyset$ and $X$ is minimal,
Theorem~\ref{thm:main} reduces to
\[
\chi(X)
=
\frac{1}{2\pi}\int_X K\,dA
+\theta(X,\infty).
\]
This formula is closely related to the classical
Jorge--Meeks formula \cite{JorgeMeeks:1983} for complete
minimal surfaces of finite total curvature in
$\mathbb{R}^3$, which states that
\[
\int_M K\,dA
=
2\pi\left(\chi(M)-\sum_{j=1}^{k}I_j\right),
\]
where $I_j$ denotes the multiplicity of the $j$-th
end. Indeed, the asymptotic description of the ends
given by Jorge and Meeks shows that the density at
infinity is precisely the sum of these multiplicities:
\[
\theta(M,\infty)=\sum_{j=1}^{k}I_j.
\]
In particular, if $M$ is properly embedded, then
each end is embedded and $I_j=1$. Consequently,
$\theta(M,\infty)=k$, and our formula agrees with
the Jorge-Meeks formula.
\end{remark}

The paper is organized as follows. In Section~\ref{sec:prelim} we recall the facts on o-minimal structures, densities and curvature of surfaces used in the paper. In Section~\ref{sec:finiteness} we prove that the total absolute curvature of a definable surface is finite (Proposition~\ref{prop:finitude_K}). In Section~\ref{sec:local} we study the geometry of a definable surface near an isolated singular point and at infinity; the key results are Propositions~\ref{prop:local} and~\ref{prop:infinity}, which compute the limits of the total geodesic curvature of the links in terms of the densities. Theorem~\ref{thm:main} is proved in Section~\ref{sec:proof} and the corollaries in Section~\ref{sec:corollaries}.

\section{Preliminaries}\label{sec:prelim}

All subsets of $\R^n$ are equipped with the metric induced by the Euclidean metric. We denote by $B_r(a)$ (resp. $\overline{B_r(a)}$, $\Sn_r(a)$) the open ball (resp. closed ball, sphere) of $\R^n$ centred at $a$ with radius $r$, and $B_r:=B_r(0)$, $\Sn_r:=\Sn_r(0)$. We write $s_m:=\HH^m(\Sn^m)$ for the volume of the unit $m$-sphere (with $s_0=2$) and $\mu_m$ for the volume of the unit $m$-ball. Recall that
\begin{equation}\label{eq:spheres}
s_{m+1}=\frac{2\pi}{m}\,s_{m-1}\qquad (m\geq1).
\end{equation}

\subsection{O-minimal structures}

\begin{definition}\label{definivel}
A \emph{structure} on $\R$ is a collection $\mathcal{S}=\{\mathcal{S}_n\}_{n\in \mathbb{N}}$, where each $\mathcal{S}_n$ is a set of subsets of $\R^n$, satisfying the following axioms:
\begin{enumerate}
\item all algebraic subsets of $\R^n$ are in $\mathcal{S}_n$;
\item for every $n$, $\mathcal{S}_n$ is a Boolean subalgebra of the power set of $\R^n$;
\item if $A\in \mathcal{S}_m$ and $B \in \mathcal{S}_n$, then $A \times B \in \mathcal{S}_{m+n}$;
\item if $\pi \colon \R^{n+1} \to \R^n$ is the projection on the first $n$ coordinates and $A\in \mathcal{S}_{n+1}$, then $\pi(A)\in \mathcal{S}_n$.
\end{enumerate}
The structure $\mathcal{S}$ is said to be \emph{o-minimal} if, moreover,
\begin{enumerate}[resume]
\item the elements of $\mathcal{S}_1$ are precisely the finite unions of points and intervals.
\end{enumerate}
An element of $\mathcal{S}_n$ is called \emph{definable} (in $\mathcal{S}$). A map $f\colon A\subset\R^n\to\R^m$ is \emph{definable} if its graph is definable.
\end{definition}

\begin{remark}\label{rem:examples}
Classical examples of o-minimal structures are: the semialgebraic sets (by the Tarski--Seidenberg theorem); the globally subanalytic sets (by Gabrielov's complement theorem; see \cite{Dries:1986}); and the structure $\R_{\exp}$ generated by the semialgebraic sets and the graph of the exponential function, whose o-minimality is due to Wilkie \cite{Wilkie:1996}.
\end{remark}

Throughout this paper we fix an o-minimal structure $\mathcal{S}$ on $\R$, and \emph{definable} means definable in $\mathcal{S}$. We shall use the following standard facts.

\begin{enumerate}[label=(O\arabic*),leftmargin=*]
\item\label{O:mono} (Monotonicity) If $f\colon(a,b)\to\R$ is definable, there are $a=a_0<a_1<\cdots<a_k=b$ such that on each $(a_{j-1},a_j)$ the function $f$ is $C^1$ and either constant or strictly monotone \cite[Ch.~3, \S1 and Ch.~7]{Dries:1998}. In particular, $\lim_{t\to a^+}f(t)$ and $\lim_{t\to b^-}f(t)$ exist in $[-\infty,+\infty]$.
\item\label{O:choice} (Definable choice and curve selection) If $A\subset\R^m\times\R^n$ is definable and $\pi(A)$ denotes its projection on $\R^m$, there is a definable map $g\colon\pi(A)\to\R^n$ with $(y,g(y))\in A$ for all $y\in\pi(A)$ \cite[Ch.~6, \S1]{Dries:1998}.
\item\label{O:dim} (Dimension) Definable sets have a well defined dimension, which coincides with the dimension as a manifold for definable $C^1$ submanifolds. If $f\colon A\to B$ is a definable map, then for each $d$ the set $B_d=\{b\in B:\dim f^{-1}(b)=d\}$ is definable and $\dim f^{-1}(B_d)=\dim B_d+d$ \cite[Ch.~4, \S1]{Dries:1998}. A definable subset of a $C^1$ manifold of dimension $m$ whose dimension is smaller than $m$ has zero $\HH^m$-measure, since it is a finite union of $C^1$ submanifolds of dimension smaller than $m$ \cite[Ch.~7]{Dries:1998}.
\item\label{O:unif} (Uniform finiteness) If $A\subset\R^m\times\R^n$ is definable, there is $N\in\mathbb{N}$ such that, for every $y\in\R^m$, the fibre $A_y=\{x\in\R^n:(y,x)\in A\}$ has at most $N$ connected components; in particular, if $A_y$ is finite, then $\#A_y\leq N$ \cite[Ch.~3, \S\S2--3]{Dries:1998}, \cite{Coste:2000}. Moreover, every definable set has finitely many connected components, and they are definable.
\item\label{O:diff} (Differentiability) The partial derivatives of a definable $C^1$ map are definable; if $N\subset\R^n$ is a definable $C^1$ submanifold, the map $x\mapsto T_xN$ (equivalently, $x\mapsto P_x$, the orthogonal projection onto $T_xN$) is definable \cite[Ch.~7]{Dries:1998}. For every $k\geq1$ and every definable $X$, the set $\Reg_k(X)$ of Definition~\ref{def:sing} is definable and $\dim\Sing_k(X)<\dim X$ \cite{DriesMiller:1996}.
\end{enumerate}

\subsection{Singular set and surfaces}

\begin{definition}\label{def:sing}
Let $X\subset \R^n$ be a subset and $k$ a positive integer. The \emph{$C^k$ singular set} of $X$, denoted by $\Sing_k(X)$, is the set of points $x\in X$ such that $U\cap X$ is not a $C^k$ submanifold of $\R^n$ for any open neighbourhood $U$ of $x$. We set $\Reg_k(X):=X\setminus\Sing_k(X)$. If $p\in\Reg_k(X)$, there is an open neighbourhood $U\subset\R^n$ of $p$ such that $X\cap U$ is a connected $C^k$ submanifold, and we set $\dim_pX:=\dim (X\cap U)$. We say that $X$ is \emph{pure $d$-dimensional} if $\dim_pX=d$ for every $p\in\Reg_1(X)$.
\end{definition}

Throughout the paper, $\Sing(X):=\Sing_2(X)$ and $\Reg(X):=\Reg_2(X)$. Note that if $X$ is pure two-dimensional, then $X$ has no isolated points (an isolated point would be a point of $\Reg_1(X)$ at which $X$ has dimension $0$), so every $p\in X$ belongs to $\overline{X\setminus\{p\}}$.

\subsection{Local density and density at infinity}\label{subsec:density}

We introduce the local density and the density at infinity; for further details we refer to \cite{KurdykaRaby:1989,SilvaSampaio:2026}.

\begin{definition}
Let $X\subset \R^n$ be a set, $d\in\mathbb{N}$ and $p\in\R^n$. We say that $X$ has \emph{$d$-dimensional density at $p$} if the limit
\[
\theta^d(X,p):=\lim_{r\to0^+}\frac{\HH^d(X\cap B_r(p))}{\mu_d r^d}
\]
exists, and that $X$ has \emph{$d$-dimensional density at infinity} if the limit
\[
\theta^d(X,\infty):=\lim_{r\to+\infty}\frac{\HH^d(X\cap B_r)}{\mu_d r^d}
\]
exists. If $X$ is a surface, we write $\theta(X,p):=\theta^2(X,p)$ and $\theta(X,\infty):=\theta^2(X,\infty)$.
\end{definition}

Since $B_{r-|a|}\subset B_r(a)\subset B_{r+|a|}$, the density at infinity does not depend on the choice of the centre of the balls. If $X$ is a $C^1$ surface near $p\in X$, then $\theta(X,p)=1$, since near $p$ it is the graph over $p+T_pX$ of a $C^1$ map with vanishing derivative at $p$.

\subsection{Curvature of surfaces}

Let $M\subset\R^n$ be a $C^2$ surface (two-dimensional submanifold). We denote by $D$ the Euclidean connection, by $\nabla$ the Levi-Civita connection of $M$, by $\nabla^\perp$ the normal connection and by $N_xM$ the normal space of $M$ at $x$. For tangent vector fields $V,W$ and normal vector fields $\nu$ we have the Gauss and Weingarten formulas (see \cite[Ch.~8]{LeeRiemannian})
\[
D_VW=\nabla_VW+\II(V,W),\qquad D_V\nu=-A_\nu V+\nabla^\perp_V\nu,
\]
where $\II$ is the second fundamental form (symmetric, bilinear and continuous since $M$ is $C^2$) and the shape operator $A_\nu\colon T_xM\to T_xM$ is the self-adjoint operator given by $\langle A_\nu v,w\rangle=\langle\II(v,w),\nu\rangle$; it depends linearly on $\nu\in N_xM$. The Gaussian curvature of $M$ is defined by the Gauss equation: if $e_1,e_2$ is an orthonormal basis of $T_xM$ and $\nu_1,\ldots,\nu_{n-2}$ is an orthonormal basis of $N_xM$, then
\begin{equation}\label{eq:gauss-eq}
K(x)=\langle\II(e_1,e_1),\II(e_2,e_2)\rangle-|\II(e_1,e_2)|^2=\sum_{\alpha=1}^{n-2}\det A_{\nu_\alpha}.
\end{equation}
If $M$ is of class $C^3$, this is the intrinsic Gaussian curvature of the induced metric, by Gauss' Theorema Egregium \cite[Ch.~8]{LeeRiemannian}.

Let $\Gamma\subset M$ be a $C^2$ curve with unit tangent vector $T$ (defined up to sign) and arc length $s$. Its \emph{curvature vector} $\kv:=\frac{dT}{ds}$ does not depend on the orientation. A \emph{unit conormal} along $\Gamma$ is a continuous unit vector field $\eta$ tangent to $M$ and orthogonal to $\Gamma$. The \emph{geodesic curvature} of $\Gamma$ with respect to $\eta$ is
\[
k_g^\eta:=\langle\kv,\eta\rangle=\langle\nabla_TT,\eta\rangle,
\]
where the second equality follows from the Gauss formula, since $\II(T,T)\perp T_xM$. Clearly $k_g^{-\eta}=-k_g^\eta$ and $|k_g^\eta|\leq|\kv|$.

\begin{theorem}[Gauss--Bonnet]\label{thm:GB-classical}
Let $S\subset\R^n$ be a compact $C^2$ surface with $C^2$ boundary, and let $\eta_{\rm in}$ be the inward pointing unit conormal along $\partial S$. Then
\[
\int_SK\dd A+\int_{\partial S}k_g^{\eta_{\rm in}}\dd s=2\pi\chi(S).
\]
\end{theorem}

For smooth oriented compact surfaces with boundary this is classical, see e.g. \cite[Ch.~9]{LeeRiemannian}. The general case reduces to this one as follows. Let $S$ be as in the theorem. By \cite[Ch.~2]{Hirsch:1976}, $S$ has a $C^\infty$ structure compatible with its $C^2$ structure, and the inclusion $\iota\colon S\to\R^n$ can be approximated in the $C^2$ topology by $C^\infty$ embeddings $\iota_k\colon S\to\R^n$. The Gaussian curvature, the curvature vector of the boundary, the inward conormal, and the area and length elements of $\iota_k(S)$ depend continuously on the $2$-jet of $\iota_k$; hence they converge uniformly to those of $S$, and it suffices to prove the formula for the smooth surfaces $\iota_k(S)$. If such a surface is not orientable, we apply the oriented case to its orientation double cover, endowed with the pulled back metric: its Euler characteristic and both integrals are twice those of the surface.

\subsection{Area formula}

We shall use the area formula in the following form (see \cite[3.2.3 and 3.2.5]{Federer:1969}, applied in charts and summed over a countable Borel partition): if $F\colon N\to N'$ is a locally Lipschitz map between $C^1$ Riemannian manifolds of the same dimension $m$ (second countable), with Jacobian $JF$, and $A\subset N$ is a Borel set, then
\begin{equation}\label{eq:area-formula}
\int_A JF\dd \HH^m=\int_{N'}\#\big(A\cap F^{-1}(y)\big)\dd\HH^m(y).
\end{equation}
On a $C^1$ submanifold of $\R^n$, $\HH^m$ coincides with the Riemannian volume \cite[\S3.2]{Federer:1969}.

\section{Finiteness of the total curvature}\label{sec:finiteness}

The finiteness of the total absolute curvature of definable submanifolds is well known; see, e.g., Fu \cite{Fu:1994} for compact subanalytic sets, Bröcker and Kuppe \cite{BrockerKuppe:2000} for compact definable sets, and Dutertre \cite{Dutertre:2004} and Grandjean \cite[Proposition 3.1]{Grandjean:2008} for (possibly noncompact) definable hypersurfaces. Since we need it for noncompact surfaces of arbitrary codimension, we include the short proof, which combines the classical argument of Chern and Lashof \cite{ChernLashof:1957} with the uniform finiteness property~\ref{O:unif}. Dutertre pointed out to us that the proofs in \cite{Dutertre:2004} and \cite{Grandjean:2008} should work for the regular part of any definable set.

We first recall an elementary identity.

\begin{lemma}\label{lem:moments}
For $q\geq1$ and $\alpha,\beta\in\{1,\ldots,q\}$,
$
\int_{\Sn^{q-1}}t_\alpha t_\beta\dd\sigma(t)=\frac{s_{q-1}}{q}\,\delta_{\alpha\beta}.
$
\end{lemma}
\begin{proof}
If $\alpha\neq\beta$, the reflection $t\mapsto(t_1,\ldots,-t_\alpha,\ldots,t_q)$ preserves $\sigma$ and changes the sign of $t_\alpha t_\beta$, so the integral vanishes. The integrals of $t_1^2,\ldots,t_q^2$ coincide by symmetry, and their sum is $\int_{\Sn^{q-1}}1\dd\sigma=s_{q-1}$.
\end{proof}

\begin{proposition}\label{prop:finitude_K}
Let $X\subset\R^n$ be a definable set and let $M$ be the union of the two-dimensional connected components of $\Reg_2(X)$. Then there is an integer $N_0$ such that
\[
\int_M |K|\dd A\leq 2\pi N_0<\infty.
\]
In particular, $K\in L^1(M)$.
\end{proposition}

\begin{proof}
$M$ is a definable $C^2$ surface by~\ref{O:diff}. If $n=2$, then $M$ is open in $\R^2$ and $K\equiv0$; so assume $n\geq3$ and put $q=n-2$.

\emph{A pointwise inequality.} Fix $x\in M$ and an orthonormal basis $\nu_1,\ldots,\nu_q$ of $N_xM$. Every unit normal vector is $\nu=\sum_\alpha t_\alpha\nu_\alpha$ with $t\in\Sn^{q-1}$, and $A_\nu=\sum_\alpha t_\alpha A_{\nu_\alpha}$. Writing the matrix of $A_{\nu_\alpha}$ in an orthonormal basis of $T_xM$ as $\begin{pmatrix}a_\alpha&b_\alpha\\ b_\alpha&c_\alpha\end{pmatrix}$, we get
\[
\det A_\nu=\Big(\sum_\alpha t_\alpha a_\alpha\Big)\Big(\sum_\beta t_\beta c_\beta\Big)-\Big(\sum_\alpha t_\alpha b_\alpha\Big)^2=\sum_{\alpha,\beta}t_\alpha t_\beta\,(a_\alpha c_\beta-b_\alpha b_\beta).
\]
Since $t\mapsto\sum_\alpha t_\alpha\nu_\alpha$ is an isometry from $\Sn^{q-1}$ onto the unit sphere $S(N_xM)$ of $N_xM$, Lemma~\ref{lem:moments} and \eqref{eq:gauss-eq} give
\[
\int_{S(N_xM)}\det A_\nu\dd\sigma_x(\nu)=\frac{s_{q-1}}{q}\sum_{\alpha}(a_\alpha c_\alpha-b_\alpha^2)=\frac{s_{q-1}}{q}\,K(x).
\]
Hence
\begin{equation}\label{eq:pointwise}
|K(x)|\leq\frac{q}{s_{q-1}}\int_{S(N_xM)}|\det A_\nu|\dd\sigma_x(\nu).
\end{equation}

\emph{The generalized Gauss map.} Let $W\subset M$ be an open set on which there is a $C^1$ orthonormal frame $\nu_1,\ldots,\nu_q$ of the normal bundle (such frames exist locally, since $NM$ is a $C^1$ vector bundle). Consider the $C^1$ map
\[
G\colon W\times\Sn^{q-1}\to\Sn^{n-1},\qquad G(x,t)=\sum_{\alpha=1}^q t_\alpha\nu_\alpha(x),
\]
and endow $W\times\Sn^{q-1}$ with the product metric. Fix $(x,t)$ and put $\nu=G(x,t)$. For $v\in T_xM$ and $w\in T_t\Sn^{q-1}=t^\perp$, the Weingarten formula gives
\[
dG_{(x,t)}(v,w)=\sum_\alpha t_\alpha D_v\nu_\alpha+\sum_\alpha w_\alpha\nu_\alpha(x)
=-A_\nu v+B(v)+\iota(w),
\]
where $B(v)=\sum_\alpha t_\alpha\nabla^\perp_v\nu_\alpha$ and $\iota(w)=\sum_\alpha w_\alpha\nu_\alpha(x)$. Differentiating $\langle\nu_\alpha,\nu_\beta\rangle=\delta_{\alpha\beta}$ we see that the matrix $\big(\langle\nabla^\perp_v\nu_\alpha,\nu_\beta\rangle\big)_{\alpha,\beta}$ is skew-symmetric, hence
\[
\langle B(v),\nu\rangle=\sum_{\alpha,\beta}t_\alpha t_\beta\langle\nabla^\perp_v\nu_\alpha,\nu_\beta\rangle=0,
\]
i.e., $B(v)\in N_xM\cap\nu^\perp$. Moreover, $\iota$ is a linear isometry from $t^\perp$ onto $N_xM\cap\nu^\perp$. With respect to the orthogonal decompositions
\[
T_{(x,t)}\big(W\times\Sn^{q-1}\big)=T_xM\oplus t^\perp
\qquad\text{and}\qquad
T_\nu\Sn^{n-1}=\nu^\perp=T_xM\oplus(N_xM\cap\nu^\perp),
\]
the matrix of $dG_{(x,t)}$ is block lower triangular,
\[
\begin{pmatrix}-A_\nu&0\\ B&\iota\end{pmatrix},
\]
so that $JG(x,t)=|\det A_\nu|\,|\det\iota|=|\det A_\nu|$. On the other hand, for $u\in\Sn^{n-1}$ we have $G(x,t)=u$ if and only if $u\perp T_xM$ and $t_\alpha=\langle u,\nu_\alpha(x)\rangle$ for all $\alpha$; hence $\#\big((Y\times\Sn^{q-1})\cap G^{-1}(u)\big)=\#\{x\in Y: u\perp T_xM\}$ for every $Y\subset W$.

Cover $M$ by countably many such open sets $W_k$ and choose a Borel partition $M=\bigsqcup_kY_k$ with $Y_k\subset W_k$. Applying \eqref{eq:area-formula} to each $Y_k\times\Sn^{q-1}$, using Fubini's theorem and summing (monotone convergence), we obtain
\begin{equation}\label{eq:chern-lashof}
\int_M\int_{S(N_xM)}|\det A_\nu|\dd\sigma_x(\nu)\dd A(x)=\int_{\Sn^{n-1}}\mathcal{N}(u)\dd\HH^{n-1}(u),
\end{equation}
where $\mathcal{N}(u):=\#\{x\in M:u\perp T_xM\}$ is the number of critical points of the height function $x\mapsto\langle x,u\rangle$ on $M$.

\emph{O-minimal bound.} The set $\mathcal{A}=\{(u,x)\in\Sn^{n-1}\times M:u\perp T_xM\}$ is definable by~\ref{O:diff}, and it is a $C^1$ manifold of dimension $2+(q-1)=n-1$ (it is diffeomorphic to the unit normal bundle of $M$). Let $E=\{u\in\Sn^{n-1}:\dim\mathcal{A}_u\geq1\}$. By~\ref{O:dim}, for each $d\geq1$ the set $E_d=\{u:\dim\mathcal{A}_u=d\}$ is definable and $\dim E_d+d=\dim(\mathcal{A}\cap(E_d\times M))\leq n-1$, so $\dim E\leq n-2$ and $\HH^{n-1}(E)=0$. For $u\notin E$ the fibre $\mathcal{A}_u$ is finite, and by~\ref{O:unif} there is $N_0$ with $\mathcal{N}(u)=\#\mathcal{A}_u\leq N_0$ for all $u\notin E$.

Combining \eqref{eq:pointwise}, \eqref{eq:chern-lashof} and \eqref{eq:spheres} (with $m=q$) we conclude
\[
\int_M|K|\dd A\leq\frac{q}{s_{q-1}}\int_{\Sn^{n-1}}\mathcal{N}\dd\HH^{n-1}\leq\frac{q\,s_{n-1}}{s_{q-1}}N_0=2\pi N_0.\qedhere
\]
\end{proof}

\begin{remark}
The proof gives the Chern--Lashof type inequality $\int_M|K|\dd A\leq\frac{2\pi}{s_{n-1}}\int_{\Sn^{n-1}}\mathcal{N}(u)\dd\HH^{n-1}(u)$, i.e., $\frac{1}{2\pi}\int_M|K|\dd A$ is bounded by the average number of critical points of the height functions on $M$.
\end{remark}

\section{Local geometry near an isolated singular point and at infinity}\label{sec:local}

\subsection{Radial tangency}

For $a\in\R^n$ we write $\rho_a(x)=|x-a|$ and $\rho=\rho_0$. If $N\subset\R^n$ is a $C^1$ submanifold, $P_x$ denotes the orthogonal projection onto $T_xN$ and $\nabla_N\rho_a(x)=P_x\big(\frac{x-a}{|x-a|}\big)$ is the gradient of $\rho_a|_N$ at $x\neq a$. Note that $|\nabla_N\rho_a(x)|=\max\{\langle\frac{x-a}{|x-a|},v\rangle: v\in T_xN,\ |v|=1\}$ is the cosine of the angle between $x-a$ and $T_xN$.

The following lemma is well known; for $a\in\overline N$ it expresses the Whitney (b)-regularity of the pair $(N,\{a\})$, which holds in o-minimal structures (see Loi \cite{Loi:1998}); see also \cite[proof of Proposition 2.1]{DinhNguyen:2023}. We include a short proof.

\begin{lemma}\label{lem:radial}
Let $N\subset\R^n$ be a definable $C^1$ submanifold of positive dimension.
\begin{enumerate}[label=\rm(\roman*)]
\item If $a\in\overline{N}$, then $|\nabla_N\rho_a(x)|\to1$ as $x\to a$, $x\in N\setminus\{a\}$.
\item If $N$ is unbounded, then $|\nabla_N\rho(x)|\to1$ as $|x|\to\infty$, $x\in N$.
\end{enumerate}
\end{lemma}

\begin{proof}
(i) Suppose the contrary. Then there is $\delta\in(0,1)$ such that $a\in\overline{Z}$, where $Z=\{x\in N\setminus\{a\}:|\nabla_N\rho_a(x)|\leq1-\delta\}$. The set $Z$ is definable by~\ref{O:diff}, so $\rho_a(Z)\subset(0,\infty)$ is a definable set with $0$ in its closure; by o-minimality it contains an interval $(0,\epsilon_1)$. By~\ref{O:choice} there is a definable map $\gamma\colon(0,\epsilon_1)\to Z$ with $|\gamma(t)-a|=t$, and by~\ref{O:mono}, after shrinking $\epsilon_1$, $\gamma$ is $C^1$ and each coordinate $\gamma_j'(t)$ has a limit $L_j\in[-\infty,+\infty]$ as $t\to0^+$. Since $\gamma_j(t)-a_j\to0$, l'H\^opital's rule gives $\lim_{t\to0^+}(\gamma_j(t)-a_j)/t=L_j$; as $|\gamma_j(t)-a_j|\leq t$, $L_j$ is finite. Thus $\gamma'(t)\to L$ and $(\gamma(t)-a)/t\to L$, with $|L|=1$. Since $\gamma'(t)\in T_{\gamma(t)}N$ and $\gamma'(t)\neq0$ for $t$ small,
\[
|\nabla_N\rho_a(\gamma(t))|\geq\Big\langle\frac{\gamma(t)-a}{t},\frac{\gamma'(t)}{|\gamma'(t)|}\Big\rangle\longrightarrow\langle L,L\rangle=1,
\]
which contradicts $\gamma(t)\in Z$.

(ii) The proof is the same: if the claim fails, the corresponding set $Z$ is unbounded, $\rho(Z)$ contains an interval $(t_1,+\infty)$, and we get a definable $C^1$ curve $\gamma\colon(t_1,\infty)\to Z$ with $|\gamma(t)|=t$ and $\gamma'(t)\to L\in[-\infty,\infty]^n$; l'H\^opital's rule (for the quotient $\gamma_j(t)/t$, whose denominator tends to $+\infty$) gives $\gamma(t)/t\to L$, so $|L|=1$, and we reach the same contradiction.
\end{proof}

\subsection{Structure of the links}

\begin{setting}\label{setting}
Let $X\subset\R^n$ be a closed definable set. We consider two situations.
\begin{enumerate}
\item[(L)] (\emph{Local case}) $p\in X$ satisfies $p\in\overline{X\setminus\{p\}}$, and there is $\varepsilon_0>0$ such that $X\cap B_{\varepsilon_0}(p)\setminus\{p\}$ is a two-dimensional $C^2$ submanifold of $\R^n$. We put $\varrho=\rho_p$ and $J=(0,\varepsilon_0)$.
\item[(I)] (\emph{Case at infinity}) $X$ is unbounded and there is $R_0>0$ such that $X\setminus\overline{B_{R_0}}$ is a two-dimensional $C^2$ submanifold of $\R^n$. We put $\varrho=\rho$ and $J=(R_0,+\infty)$.
\end{enumerate}
In both cases, for an interval $J'\subset J$ we write $U_{J'}=X\cap\varrho^{-1}(J')$, which is an open subset of a $C^2$ surface; for $t\in J$ we write
\[
L_t=X\cap\varrho^{-1}(t)\qquad\text{and}\qquad\ell(t)=\HH^1(L_t).
\]
We write $\nabla\varrho$ for the gradient of $\varrho|_{U_J}$, $u=\nabla\varrho/|\nabla\varrho|$ wherever $\nabla\varrho\neq0$, and $V=\nabla\varrho/|\nabla\varrho|^2$. Since $U_J$ is $C^2$, $x\mapsto P_x$ is $C^1$, so $\nabla\varrho$, $u$ and $V$ are $C^1$ on their domains.
\end{setting}

The vector field $u$ is the unit conormal along the curves $L_t$ pointing towards increasing values of $\varrho$, i.e., away from $p$ in case (L) and away from the origin in case (I).

\begin{lemma}\label{lem:structure}
In Setting~\ref{setting}, there is an interval $J_1\subset J$, of the form $J_1=(0,\varepsilon)$ in case (L) and $J_1=(R_1,+\infty)$ in case (I), such that:
\begin{enumerate}[label=\rm(\alph*)]
\item $|\nabla\varrho|\geq\frac12$ on $U_{J_1}$, and $\delta(t):=\sup_{x\in L_t}\big||\nabla\varrho(x)|^{-1}-1\big|\to0$ as $t\to0^+$ in case (L), and as $t\to+\infty$ in case (I);
\item for every $t\in J_1$, $L_t$ is a nonempty compact $C^2$ curve, i.e., a finite disjoint union of $C^2$ circles;
\item for every $t_0\in J_1$, there is a $C^1$ diffeomorphism $h\colon J_1\times L_{t_0}\to U_{J_1}$ such that $\varrho(h(t,q))=t$, $h(t_0,q)=q$ and $\partial_th=V\circ h$; in particular, each $h(t,\cdot)\colon L_{t_0}\to L_t$ is a $C^1$ diffeomorphism, all the curves $L_t$, $t\in J_1$, have the same number $b$ of connected components, and $U_{(a,a')}$ is diffeomorphic to $(a,a')\times L_{t_0}$ for $(a,a')\subset J_1$;
\item in case (L), for every $r\in(0,\varepsilon]$, $X\cap B_r(p)$ is homeomorphic to the open cone $\Cone(L_{t_0})=\big([0,r)\times L_{t_0}\big)/\big(\{0\}\times L_{t_0}\big)$; in particular, it is contractible.
\end{enumerate}
\end{lemma}

\begin{proof}
(a) Apply Lemma~\ref{lem:radial} to the definable $C^1$ submanifold $N=U_J$: in case (L), $p\in\overline{N}$ because $p\in\overline{X\setminus\{p\}}$; in case (I), $N$ is unbounded because $X$ is.

(b) By (a), $\varrho|_{U_{J_1}}$ is a $C^2$ submersion, so each $L_t$, $t\in J_1$, is a $C^2$ curve. It is compact because $X$ is closed and $L_t$ is bounded. The map $\varrho|_{U_{J_1}}\colon U_{J_1}\to J_1$ is open (it is a submersion) and closed (if $\varrho(x_k)\to t\in J_1$, then the $x_k$ eventually lie in the compact set $X\cap\varrho^{-1}([t-\tau,t+\tau])\subset U_{J_1}$ for some small $\tau>0$, so a subsequence converges to some $x\in X$ with $\varrho(x)=t$); since $U_{J_1}\neq\emptyset$ and $J_1$ is connected, $\varrho(U_{J_1})=J_1$. Finally, a compact one-dimensional manifold is a finite disjoint union of circles (see, e.g., \cite[Appendix]{Milnor:1965}).

(c) Let $\Phi$ be the maximal flow of the $C^1$ vector field $V$ on $U_{J_1}$; it is of class $C^1$ \cite[Ch.~V]{Hartman:1982}. Since $d\varrho(V)=\langle\nabla\varrho,V\rangle=1$, we have $\varrho(\Phi_s(x))=\varrho(x)+s$ while $\Phi_s(x)$ is defined (see the figure \ref{fig:radial-flow}). We claim that $\Phi_s(x)$ is defined for every $s$ with $\varrho(x)+s\in J_1$. Indeed, if the maximal interval of definition were $(\alpha,\omega)$ with $\varrho(x)+\omega\in J_1$, then $\Phi_s(x)$, $s\in[0,\omega)$, would stay in the compact set $X\cap\varrho^{-1}([\varrho(x),\varrho(x)+\omega])\subset U_{J_1}$, on which $V$ is bounded by (a); hence it could be extended beyond $\omega$ \cite[Ch.~II, \S3]{Hartman:1982}, a contradiction; the same argument applies to $\alpha$. Then $h(t,q)=\Phi_{t-t_0}(q)$ has the required properties, and its inverse $x\mapsto\big(\varrho(x),\Phi_{t_0-\varrho(x)}(x)\big)$ is $C^1$.

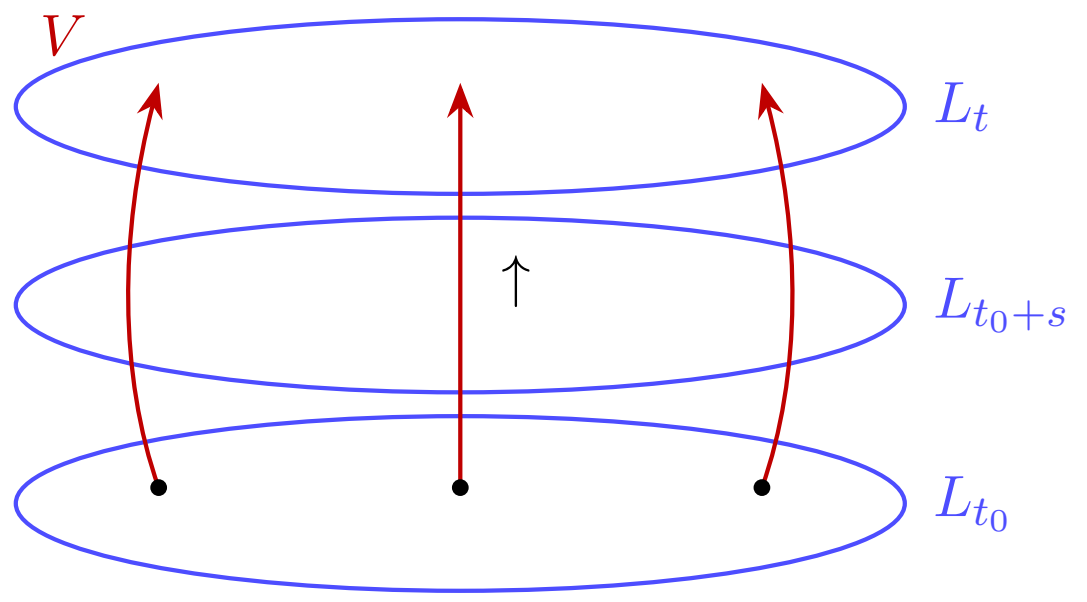
\begin{figure}[ht]
\centering
\begin{tikzpicture}[scale=1.05, >=Stealth]

    \draw[thick, blue!70]
        (0,0) ellipse (2.8 and 0.55);

    \draw[thick, blue!70]
        (0,1.25) ellipse (2.8 and 0.55);

    \draw[thick, blue!70]
        (0,2.5) ellipse (2.8 and 0.55);

    \draw[thick, red!75!black, ->]
        (-1.9,0.1) .. controls (-2.15,0.8)
        and (-2.15,1.7) .. (-1.9,2.65);

    \draw[thick, red!75!black, ->]
        (0,0.1) -- (0,2.65);

    \draw[thick, red!75!black, ->]
        (1.9,0.1) .. controls (2.15,0.8)
        and (2.15,1.7) .. (1.9,2.65);

    \fill (-1.9,0.1) circle (1.5pt);
    \fill (0,0.1) circle (1.5pt);
    \fill (1.9,0.1) circle (1.5pt);

    \node[blue!70, right] at (2.85,0)
        {$L_{t_0}$};

    \node[blue!70, right] at (2.85,1.25)
        {$L_{t_0+s}$};

    \node[blue!70, right] at (2.85,2.5)
        {$L_t$};

    \node[red!75!black] at (-2.5,2.95)
        {$V$};

    \node at (0.35,1.4) {$\uparrow$};

\end{tikzpicture}
\caption{Schematic representation of the flow of $V$,
which transports points between the level curves of
$\varrho$.}
\label{fig:radial-flow}
\end{figure}

(d) Define $\hat h\colon[0,r)\times L_{t_0}\to X\cap B_r(p)$ by $\hat h(0,q)=p$ and $\hat h=h$ on $(0,r)\times L_{t_0}$. Since $|\hat h(t,q)-p|=t$, $\hat h$ is continuous, and it induces a continuous bijection $\Cone(L_{t_0})\to X\cap B_r(p)$. For every $r'<r$, its restriction to the compact image of $[0,r']\times L_{t_0}$ is a homeomorphism onto $X\cap\overline{B_{r'}(p)}$ (a continuous bijection from a compact space onto a Hausdorff space); hence the bijection is a homeomorphism.
\end{proof}

\begin{lemma}\label{lem:total-curvature}
In Setting~\ref{setting}, let $J_1$ be as in Lemma~\ref{lem:structure}, and for $x\in U_{J_1}$ let $\kv(x)$ be the curvature vector at $x$ of the curve $L_{\varrho(x)}$. Then $\kv\colon U_{J_1}\to\R^n$ is continuous and definable, and there is a constant $C>0$ such that
\[
\int_{L_t}|\kv|\dd s\leq C\qquad\text{for every }t\in J_1.
\]
In particular, $\int_{L_t}|k_g^u|\dd s\leq C$ for every $t\in J_1$. Moreover, in case (L), $\ell(t)\leq C_1t$ for every $t\in J_1$, where $C_1=4(b+3C/\pi)$.
\end{lemma}

\begin{proof}
For $x\in U_{J_1}$, the tangent line of $L_{\varrho(x)}$ at $x$ is $T_xU_{J_1}\cap\nabla\varrho(x)^\perp$, and the orthogonal projection onto it is $Q_x=P_x-u(x)\otimes u(x)$, where $(a\otimes a)(w)=\langle a,w\rangle a$. Thus $x\mapsto Q_x$ is $C^1$ and definable. Near each point there is a $C^1$ unit vector field $T$ spanning these tangent lines, unique up to sign, and $\kv=dT(T)$. Since $Q=T\otimes T$, we have $dQ(v)=dT(v)\otimes T+T\otimes dT(v)$, and using $\langle dT(v),T\rangle=0$ we get
\[
\kv(x)=dQ_x(T(x))\,T(x).
\]
The right-hand side is an even function of $T(x)$; hence $\kv$ is a well defined continuous map, and it is definable, its graph being $\{(x,w):\exists\, v\in\operatorname{Im}Q_x,\ |v|=1,\ w=dQ_x(v)v\}$, where $dQ$ is definable by~\ref{O:diff}.

Fix $j\in\{1,\ldots,n\}$ and let $Z_j=\{x\in U_{J_1}:\kappa_j(x)=0\}$, where $\kappa_j$ is the $j$-th coordinate of $\kv$. The family $\{L_t\setminus Z_j\}_{t\in J_1}$ is definable, so by~\ref{O:unif} there is $N_j$ such that $L_t\setminus Z_j$ has at most $N_j$ connected components for every $t$. Let $I$ be one of them; it is a connected open subset of some component $\Sigma$ of $L_t$. Let $\sigma$ be an arc length parametrization of $\Sigma$, so that $\sigma''=\kv\circ\sigma$. On $\sigma^{-1}(I)$ the continuous function $\sigma_j''=\kappa_j\circ\sigma$ does not vanish, so it has constant sign. If $I=\Sigma$, then $\int_\Sigma\kappa_j\dd s=\int\sigma_j''=0$ by periodicity, which is impossible. Hence $I=\sigma((c,c'))$ is an open arc and
\[
\int_I|\kappa_j|\dd s=\Big|\int_c^{c'}\sigma_j''\Big|=|\sigma_j'(c')-\sigma_j'(c)|\leq2.
\]
Since $\kappa_j=0$ on $Z_j$, we get $\int_{L_t}|\kappa_j|\dd s\leq2N_j$, and therefore $\int_{L_t}|\kv|\dd s\leq\sum_j\int_{L_t}|\kappa_j|\dd s\leq C:=2\sum_jN_j$. The second assertion follows from $|k_g^u|\leq|\kv|$.

Assume now we are in case (L), and let $\Sigma$ be a component of $L_t$, of length $\lambda$ and total curvature $c_\Sigma=\int_\Sigma|\kv|\dd s$, and $\sigma\colon[0,\lambda]\to\Sigma$ an arc length parametrization. Choose $0=\tau_0<\tau_1<\cdots<\tau_m=\lambda$ with $m\leq1+3c_\Sigma/\pi$ and $\int_{\tau_{k-1}}^{\tau_k}|\sigma''|\leq\pi/3$ for every $k$ (this is possible since $\tau\mapsto\int_0^\tau|\sigma''|$ is continuous and nondecreasing). For $\tau\in[\tau_{k-1},\tau_k]$, the spherical distance between $\sigma'(\tau)$ and $\sigma'(\tau_{k-1})$ is at most the length $\int_{\tau_{k-1}}^\tau|\sigma''|\leq\pi/3$ of the curve $\sigma'$ on $\Sn^{n-1}$, so $\langle\sigma'(\tau),\sigma'(\tau_{k-1})\rangle\geq\frac12$. Since $\Sigma\subset\Sn_t(p)$ has diameter at most $2t$,
\[
2t\geq\langle\sigma(\tau_k)-\sigma(\tau_{k-1}),\sigma'(\tau_{k-1})\rangle=\int_{\tau_{k-1}}^{\tau_k}\langle\sigma'(\tau),\sigma'(\tau_{k-1})\rangle\dd\tau\geq\frac{\tau_k-\tau_{k-1}}{2}.
\]
Thus $\lambda\leq4tm\leq4t(1+3c_\Sigma/\pi)$. Summing over the $b$ components of $L_t$ we get $\ell(t)\leq4t(b+3C/\pi)$.
\end{proof}

\begin{lemma}\label{lem:variation}
In Setting~\ref{setting}, let $J_1$, $\delta$ and $C$ be as in Lemmas~\ref{lem:structure} and~\ref{lem:total-curvature}, and put
\[
I(t):=\int_{L_t}k_g^u\dd s=\int_{L_t}\langle\kv,u\rangle\dd s,\qquad t\in J_1.
\]
Then:
\begin{enumerate}[label=\rm(\roman*)]
\item (First variation) $\ell$ is of class $C^1$ on $J_1$ and $\ell'(t)=-\int_{L_t}\frac{\langle\kv,u\rangle}{|\nabla\varrho|}\dd s$; in particular, $|\ell'(t)+I(t)|\leq C\,\delta(t)$.
\item (Coarea) For $a<a'$ in $J_1$, $\HH^2(U_{(a,a')})=\int_a^{a'}\Big(\int_{L_t}\frac{\dd s}{|\nabla\varrho|}\Big)\dd t$; in particular, $\big|\HH^2(U_{(a,a')})-\int_a^{a'}\ell(t)\dd t\big|\leq\int_a^{a'}\delta(t)\ell(t)\dd t$.
\item (Gauss--Bonnet on annuli) For $a<a'$ in $J_1$, $I(a')-I(a)=\int_{U_{(a,a')}}K\dd A$.
\end{enumerate}
\end{lemma}

\begin{proof}
Fix $t_0\in J_1$ and let $h$ be as in Lemma~\ref{lem:structure}(c). Let $\Sigma_1,\ldots,\Sigma_b$ be the components of $L_{t_0}$, let $c_i\colon\R/\lambda_i\mathbb{Z}\to\Sigma_i$ be arc length parametrizations, and put $H_i(t,\vartheta)=h(t,c_i(\vartheta))$. Then $\ell(t)=\sum_i\int_0^{\lambda_i}|\partial_\vartheta H_i(t,\vartheta)|\dd\vartheta$, and $\partial_\vartheta H_i\neq0$.

(i) The map $H_i$ is $C^1$ and $\partial_tH_i=V\circ H_i$ is $C^1$; hence $\partial_\vartheta\partial_tH_i$ exists and is continuous, and by Schwarz's theorem \cite[Theorem 9.41]{Rudin:1976} also $\partial_t\partial_\vartheta H_i$ exists and equals $\partial_\vartheta\partial_tH_i$. Differentiating under the integral sign,
\[
\ell'(t)=\sum_i\int_0^{\lambda_i}\big\langle\partial_\vartheta(V\circ H_i),T\circ H_i\big\rangle\dd\vartheta,
\]
where $T\circ H_i=\partial_\vartheta H_i/|\partial_\vartheta H_i|$, and $\ell'$ is continuous. For fixed $t$, $\vartheta\mapsto T(H_i(t,\vartheta))$ is $C^1$, because $T$ is a $C^1$ unit tangent field along the $C^2$ curve $L_t$ and $H_i(t,\cdot)$ is a $C^1$ regular parametrization of a component of $L_t$; moreover $\partial_\vartheta(T\circ H_i)=|\partial_\vartheta H_i|\,\kv\circ H_i$. Integrating by parts on the circle $\R/\lambda_i\mathbb{Z}$,
\[
\ell'(t)=-\sum_i\int_0^{\lambda_i}\langle V\circ H_i,\kv\circ H_i\rangle|\partial_\vartheta H_i|\dd\vartheta=-\int_{L_t}\langle V,\kv\rangle\dd s=-\int_{L_t}\frac{\langle u,\kv\rangle}{|\nabla\varrho|}\dd s,
\]
since $V=u/|\nabla\varrho|$. Consequently
\[
|\ell'(t)+I(t)|=\Big|\int_{L_t}\langle u,\kv\rangle\big(1-|\nabla\varrho|^{-1}\big)\dd s\Big|\leq\delta(t)\int_{L_t}|\kv|\dd s\leq C\,\delta(t).
\]

(ii) Endow $(a,a')\times L_{t_0}$ with the product metric. At $(t,q)$ we have $dh(\partial_t)=V$, which is orthogonal to $dh(T_qL_{t_0})=T_{h(t,q)}L_t$; hence $Jh(t,q)=|V(h(t,q))|\,J\big(h(t,\cdot)\big)(q)$. By \eqref{eq:area-formula} applied to $h$, Fubini's theorem, and \eqref{eq:area-formula} applied to $h(t,\cdot)\colon L_{t_0}\to L_t$,
\[
\HH^2(U_{(a,a')})=\int_a^{a'}\!\!\int_{L_{t_0}}|V(h(t,q))|\,J\big(h(t,\cdot)\big)(q)\dd\HH^1(q)\dd t=\int_a^{a'}\!\!\int_{L_t}|V|\dd s\dd t,
\]
and $|V|=|\nabla\varrho|^{-1}$. The last assertion follows from $\big||\nabla\varrho|^{-1}-1\big|\leq\delta(t)$ on $L_t$.

(iii) The set $S=X\cap\varrho^{-1}([a,a'])$ is a compact $C^2$ surface with boundary $L_a\sqcup L_{a'}$ (since $\varrho|_{U_{J_1}}$ is a $C^2$ submersion), which is diffeomorphic to $[a,a']\times L_{t_0}$ by Lemma~\ref{lem:structure}(c), i.e., to a disjoint union of $b$ annuli. Hence $\chi(S)=0$. The inward unit conormal of $S$ is $u$ along $L_a$ and $-u$ along $L_{a'}$. By Theorem~\ref{thm:GB-classical},
\[
0=\int_SK\dd A+\int_{L_a}k_g^u\dd s+\int_{L_{a'}}k_g^{-u}\dd s=\int_{U_{(a,a')}}K\dd A+I(a)-I(a'),
\]
where we used that $L_a\cup L_{a'}$ has zero area.
\end{proof}

\subsection{Total geodesic curvature of the links and densities}

The following proposition is motivated by Grieser's approach in \cite{Grieser:2003}: we analyse the geodesic curvature of the links, but the asymptotic expansions coming from resolution of singularities are replaced by o-minimal finiteness arguments.

\begin{proposition}\label{prop:local}
In Setting~\ref{setting}, case (L), let $L_r=X\cap\Sn_r(p)$ and let $u$ be the unit conormal along $L_r$ pointing away from $p$. Then:
\begin{enumerate}[label=\rm(\alph*)]
\item the limit $\Theta:=-\frac{1}{2\pi}\lim_{r\to0^+}\int_{L_r}k_g^u\dd s$ exists;
\item $\lim_{r\to0^+}\frac{\HH^1(L_r)}{2\pi r}=\Theta$; in particular, $\Theta\geq0$;
\item the density $\theta(X,p)$ exists and $\theta(X,p)=\Theta$; thus
\[
\lim_{r\to0^+}\int_{L_r}k_g^u\dd s=-2\pi\,\theta(X,p);
\]
\item for every $r\in J_1$, $\displaystyle\int_{L_r}k_g^u\dd s=-2\pi\,\theta(X,p)+\int_{X\cap B_r(p)}K\dd A$, and $\HH^2(X\cap B_r(p))\leq C_1r^2$.
\end{enumerate}
\end{proposition}

\begin{proof}
Let $J_1=(0,\varepsilon)$ be as in Lemma~\ref{lem:structure}. By Proposition~\ref{prop:finitude_K} (applied to $X$), $K$ is integrable on $U_{J_1}$.

By Lemma~\ref{lem:variation}(iii), for $0<a<a'<\varepsilon$ we have $|I(a')-I(a)|\leq\int_{U_{(a,a')}}|K|\dd A$, which tends to $0$ as $a,a'\to0^+$ by the integrability of $K$. By the Cauchy criterion, $I_0:=\lim_{a\to0^+}I(a)$ exists, which proves (a) with $\Theta=-I_0/(2\pi)$. Letting $a\to0^+$ in Lemma~\ref{lem:variation}(iii) we also get
\begin{equation}\label{eq:I-local}
I(r)=I_0+\int_{U_{(0,r)}}K\dd A\qquad(r\in J_1).
\end{equation}

By Lemma~\ref{lem:variation}(i) and Lemma~\ref{lem:structure}(a), $|\ell'(t)+I(t)|\leq C\delta(t)\to0$, hence $\lim_{t\to0^+}\ell'(t)=-I_0=2\pi\Theta$.

By Lemma~\ref{lem:total-curvature}, $0\leq\ell(t)\leq C_1t$, so $\ell(t)\to0$ as $t\to0^+$. Since $\ell'$ is bounded near $0$ (it has a limit), we get $\ell(t)=\int_0^t\ell'(\tau)\dd\tau$, and therefore $\ell(t)/t\to2\pi\Theta$. This proves (b).

Since $X\cap B_r(p)=\{p\}\cup U_{(0,r)}$, letting $a\to0^+$ in Lemma~\ref{lem:variation}(ii) (monotone convergence) gives
\[
\Big|\HH^2(X\cap B_r(p))-\int_0^r\ell(t)\dd t\Big|\leq\int_0^r\delta(t)\ell(t)\dd t\leq\Big(\sup_{(0,r)}\delta\Big)\frac{C_1r^2}{2}=o(r^2).
\]
By (b), $\int_0^r\ell(t)\dd t=\pi\Theta r^2+o(r^2)$. Hence $\HH^2(X\cap B_r(p))/(\pi r^2)\to\Theta$, which proves (c). Finally, (d) follows from \eqref{eq:I-local} and $I_0=-2\pi\theta(X,p)$, together with $\HH^2(X\cap B_r(p))\leq\int_0^r2\ell(t)\dd t\leq C_1r^2$ (note $|\nabla\varrho|^{-1}\leq2$).
\end{proof}

\begin{remark}\label{rem:local}
(1) Applying Proposition~\ref{prop:local} to the closure of each connected component $X^j$ of $U_{J_1}$, we see that the density of $X$ at $p$ is the sum of the densities of its local ``branches'' $\overline{X^j}$, and that $\lim_{r\to0^+}\int_{L_r\cap X^j}k_g^u\dd s=-2\pi\,\theta(\overline{X^j},p)$.

(2) Part (b) is a special case of \cite[Proposition 2.1]{DinhNguyen:2023}, where the density of a definable set of any dimension is expressed by the volumes of its links. The existence of the density of subanalytic sets is due to Kurdyka and Raby \cite{KurdykaRaby:1989} (for complex analytic sets, to Thie \cite{Thie:1967}); for sets definable in an arbitrary
o-minimal structure, we refer to Sampaio and Silva
\cite{SilvaSampaio:2026} and Dinh and Nguyen
\cite{DinhNguyen:2023}.
\end{remark}

The behaviour at infinity is analysed in exactly the same way; the only difference is the sign of the conormal that points into the compact part.

\begin{proposition}\label{prop:infinity}
In Setting~\ref{setting}, case (I), let $L_R=X\cap\Sn_R$ and let $-u$ be the unit conormal along $L_R$ pointing towards the origin. Then:
\begin{enumerate}[label=\rm(\alph*)]
\item the limit $\Theta_\infty:=\frac{1}{2\pi}\lim_{R\to+\infty}\int_{L_R}k_g^{-u}\dd s$ exists;
\item $\lim_{R\to+\infty}\frac{\HH^1(L_R)}{2\pi R}=\Theta_\infty$; in particular, $\Theta_\infty\geq0$;
\item if $\HH^2(X\cap B_{R})<\infty$ for some $R\in J_1$, then the density at infinity $\theta(X,\infty)$ exists and equals $\Theta_\infty$; thus
\[
\lim_{R\to+\infty}\int_{L_R}k_g^{-u}\dd s=2\pi\,\theta(X,\infty);
\]
\item for every $R\in J_1$, $\displaystyle\int_{L_R}k_g^{-u}\dd s=2\pi\,\Theta_\infty+\int_{X\setminus\overline{B_R}}K\dd A$.
\end{enumerate}
\end{proposition}

\begin{proof}
Let $J_1=(R_1,+\infty)$ be as in Lemma~\ref{lem:structure}; $K$ is integrable on $U_{J_1}$ by Proposition~\ref{prop:finitude_K}. Note that $\int_{L_R}k_g^{-u}\dd s=-I(R)$.

As in the local case, Lemma~\ref{lem:variation}(iii) and the integrability of $K$ show that $I(R)$ satisfies the Cauchy criterion as $R\to+\infty$; let $I_\infty$ be its limit and $\Theta_\infty=-I_\infty/(2\pi)$. Letting $a'\to+\infty$ in Lemma~\ref{lem:variation}(iii) we get $I_\infty-I(R)=\int_{U_{(R,\infty)}}K\dd A$, which gives (a) and (d).

By Lemma~\ref{lem:variation}(i), $\ell'(R)\to-I_\infty=2\pi\Theta_\infty$ as $R\to+\infty$. Fixing $R_2\in J_1$, from $\ell(R)=\ell(R_2)+\int_{R_2}^R\ell'(\tau)\dd\tau$ we obtain $\ell(R)/R\to2\pi\Theta_\infty$, which proves (b). In particular $\ell(t)=O(t)$ on $[R_2,\infty)$.

Assume that $\HH^2(X\cap B_{R})<\infty$ for some $R\in J_1$; by Lemma~\ref{lem:variation}(ii) this then holds for every $R\in J_1$. For $R>R_2$, since $\HH^2(L_{R_2})=0$,
\[
\HH^2(X\cap B_R)=\HH^2(X\cap B_{R_2})+\HH^2(U_{(R_2,R)}),
\]
and by Lemma~\ref{lem:variation}(ii), (b) and $\delta(t)\to0$,
\[
\HH^2(U_{(R_2,R)})=\int_{R_2}^R\ell(t)\dd t+O\Big(\int_{R_2}^R\delta(t)\,t\dd t\Big)=\pi\Theta_\infty R^2+o(R^2).
\]
Hence $\HH^2(X\cap B_R)/(\pi R^2)\to\Theta_\infty$, which proves (c).
\end{proof}

\begin{remark}\label{rem:inversion}
Proposition~\ref{prop:infinity} can also be deduced from Proposition~\ref{prop:local} via the inversion $\iota(x)=x/|x|^2$: the set $\widetilde X=\overline{\iota(X\setminus\{0\})}$ satisfies Setting~\ref{setting}(L) at $0$, $\iota$ maps $L_R$ onto $\widetilde X\cap\Sn_{1/R}$ by a homothety of ratio $R^{-2}$, so that $\HH^1(\widetilde X\cap\Sn_{1/R})/(2\pi R^{-1})=\HH^1(L_R)/(2\pi R)$ and therefore $\theta(\widetilde X,0)=\theta(X,\infty)$ by Proposition~\ref{prop:local}(b) and Proposition~\ref{prop:infinity}(b),(c); and, since $\iota$ is conformal with factor $|x|^{-2}$, the transformation rule of the geodesic curvature under conformal changes of metric gives
\[
\int_{L_R}k_g^{-u}\dd s=\int_{\widetilde X\cap\Sn_{1/R}}k_g^{\tilde u}\dd\tilde s+\frac{2}{R}\int_{L_R}|\nabla\rho|\dd s,
\]
where $\tilde u$ is the conormal pointing away from $0$. The direct proof above avoids these computations.
\end{remark}

\section{Proof of Theorem~\ref{thm:main}}\label{sec:proof}

We divide the proof into four steps. First, we choose suitable radii. Next, we compute the Euler characteristic of $X_{r,R}$. We then apply the Gauss–Bonnet theorem to $X_{r,R}$. Finally, we combine the previous results to obtain the desired formula.

Let $X$ be as in Theorem~\ref{thm:main} and $M=\Reg(X)$, which is a two-dimensional $C^2$ surface since $X$ is pure two-dimensional. By Proposition~\ref{prop:finitude_K}, $K\in L^1(M)$.

\emph{Choice of the radii.} Let $\varepsilon_0\in(0,1]$ be smaller than $\frac12|p_i-p_j|$ for all $i\neq j$. Since $X$ has no isolated points, each $p_i$ satisfies Setting~\ref{setting}(L) with this $\varepsilon_0$. Let $\varepsilon>0$ be such that the conclusions of Lemmas~\ref{lem:structure}, \ref{lem:total-curvature} and~\ref{lem:variation} and of Proposition~\ref{prop:local} hold on $(0,\varepsilon)$ for every $p_i$. If $X$ is unbounded, $X$ satisfies Setting~\ref{setting}(I) with $R_0=1+\max_i|p_i|$, and we let $R_1$ be such that the same results and Proposition~\ref{prop:infinity} hold on $(R_1,+\infty)$. Fix $r\in(0,\varepsilon)$ and $R>R_1+\varepsilon$; if $X$ is bounded, fix instead $R$ with $X\subset B_{R/2}$. Put
\[
X_{r,R}=\{x\in X:|x-p_i|\geq r\ \text{for all } i,\ |x|\leq R\}.
\]
Then $X_{r,R}\subset M$, and $X_{r,R}$ is a compact $C^2$ surface with boundary
\[
\partial X_{r,R}=\bigsqcup_{i=1}^sL_r(p_i)\ \sqcup\ L_R,\qquad L_r(p_i)=X\cap\Sn_r(p_i),\quad L_R=X\cap\Sn_R,
\]
because the corresponding distance functions are $C^2$ submersions near these curves (Lemma~\ref{lem:structure}(a)); if $X$ is bounded, $L_R=\emptyset$.

\emph{The Euler characteristic.} Consider the open subsets of $X$
\begin{align*}
\Omega_1&=\{x\in X:|x-p_i|>r/2\ \text{for all } i,\ |x|<2R\},\\
\Omega_2&=X\cap\Big(\bigcup_{i=1}^sB_r(p_i)\cup\{|x|>R\}\Big).
\end{align*}
Then $X=\Omega_1\cup\Omega_2$ and, by Lemma~\ref{lem:structure}(c),
\[
\Omega_1\cap\Omega_2=\bigsqcup_{i=1}^s\big(X\cap\{r/2<\rho_{p_i}<r\}\big)\ \sqcup\ \big(X\cap\{R<\rho<2R\}\big)
\]
is diffeomorphic to a disjoint union of products of open intervals with finite unions of circles, so $\chi(\Omega_1\cap\Omega_2)=0$. By Lemma~\ref{lem:structure}(d), each $X\cap B_r(p_i)$ is contractible, and by Lemma~\ref{lem:structure}(c), $X\cap\{|x|>R\}$ is diffeomorphic to $(R,\infty)\times L_{R}$; hence $\chi(\Omega_2)=s$. Using the product structures given by Lemma~\ref{lem:structure}(c) on the collars $X\cap\{r/2<\rho_{p_i}\leq r\}$ and $X\cap\{R\leq\rho<2R\}$, we see that $\Omega_1$ deformation retracts onto $X_{r,R}$, so $\chi(\Omega_1)=\chi(X_{r,R})$. All these spaces have finitely generated homology: $X_{r,R}$ is a compact definable set, hence definably triangulable \cite[Ch.~8]{Dries:1998}, and the others are homotopy equivalent to points or to finite unions of circles. By the Mayer--Vietoris sequence for the open cover $\{\Omega_1,\Omega_2\}$ (see \cite[\S2.2]{Hatcher:2002}), the homology of $X$ is finitely generated and
\begin{equation}\label{eq:chi-Xr}
\chi(X)=\chi(\Omega_1)+\chi(\Omega_2)-\chi(\Omega_1\cap\Omega_2)=\chi(X_{r,R})+s.
\end{equation}

\emph{ Gauss--Bonnet on $X_{r,R}$.} By Theorem~\ref{thm:GB-classical},
\begin{equation}\label{eq:gauss-bonnet-Xr}
2\pi\chi(X_{r,R})=\int_{X_{r,R}}K\dd A+\sum_{i=1}^s\int_{L_r(p_i)}k_g^{\eta_{\rm in}}\dd s+\int_{L_R}k_g^{\eta_{\rm in}}\dd s.
\end{equation}
Along $L_r(p_i)$ the inward unit conormal of $X_{r,R}$ points away from $p_i$, and along $L_R$ it points towards the origin. Hence Propositions~\ref{prop:local}(d) and~\ref{prop:infinity}(d) give
\begin{align*}
\int_{L_r(p_i)}k_g^{\eta_{\rm in}}\dd s&=-2\pi\,\theta(X,p_i)+\int_{X\cap B_r(p_i)}K\dd A\qquad(i=1,\ldots,s),\\
\int_{L_R}k_g^{\eta_{\rm in}}\dd s&=2\pi\,\Theta_\infty+\int_{X\setminus\overline{B_R}}K\dd A.
\end{align*}
Moreover, $\HH^2(X\cap B_R)<\infty$: indeed, $X\cap B_R$ is the union of the sets $X\cap B_r(p_i)$, which have finite area by Proposition~\ref{prop:local}(d), and of the compact $C^2$ surface $X_{r,R}$. Hence, by Proposition~\ref{prop:infinity}(c), $\theta(X,\infty)$ exists and equals $\Theta_\infty$. If $X$ is bounded, $\theta(X,\infty)=0$ and the term corresponding to $L_R$ is absent.

\emph{Conclusion.} Since $M$ is, up to the null set $\bigcup_iL_r(p_i)\cup L_R$, the disjoint union of $X_{r,R}\cap M$, of the sets $X\cap B_r(p_i)\setminus\{p_i\}$ and of $X\setminus\overline{B_R}$, substituting into \eqref{eq:gauss-bonnet-Xr} we get
\[
2\pi\chi(X_{r,R})=\int_XK\dd A-2\pi\sum_{i=1}^s\theta(X,p_i)+2\pi\,\theta(X,\infty).
\]
Together with \eqref{eq:chi-Xr} this gives
\[
\chi(X)
=
s+\frac{1}{2\pi}\int_X K\dd A
-\sum_{i=1}^s\theta(X,p_i)
+\theta(X,\infty),
\]
which completes the proof of Theorem~\ref{thm:main}.\qed

\begin{remark}
The proof shows slightly more: for every $r$ and $R$ as in Step~1, $\chi(X_{r,R})=\frac{1}{2\pi}\int_XK\dd A-\sum_i\theta(X,p_i)+\theta(X,\infty)$; no limit $r\to0^+$, $R\to+\infty$ is needed, since the curvature of the removed pieces is accounted for exactly by Propositions~\ref{prop:local}(d) and~\ref{prop:infinity}(d).
\end{remark}

\section{Proofs of the corollaries}\label{sec:corollaries}

\begin{proof}[Proof of Corollary~\ref{cor:minimal}]
Let $M=\Reg(X)$ and $x\in M$, and let $e_1,e_2$ be an orthonormal basis of $T_xM$. Minimality means $\II(e_1,e_1)+\II(e_2,e_2)=0$, so by \eqref{eq:gauss-eq}
\[
K(x)=-|\II(e_1,e_1)|^2-|\II(e_1,e_2)|^2=-\tfrac12|\II_x|^2\leq0.
\]
Hence $\int_XK\dd A\leq0$, and the inequality follows from Theorem~\ref{thm:main}.

If equality holds, then $\II\equiv0$ on $M$. Thus the Gauss map $x\mapsto T_xM$ is locally constant (its derivative is given by $\II$), and each connected component $M_\alpha$ of $M$ (there are finitely many by~\ref{O:unif}) is an open subset of an affine $2$-plane $P_\alpha$. Since $X$ is closed and the components of $M$ are open and closed in $M$, the set $F_\alpha=\overline{M_\alpha}\setminus M_\alpha$ is contained in the finite set $\Sing(X)$, hence $F_\alpha\subset P_\alpha$ is finite. Now $M_\alpha$ is open and closed in the connected set $P_\alpha\setminus F_\alpha$, so $M_\alpha=P_\alpha\setminus F_\alpha$ and $\overline{M_\alpha}=P_\alpha$. Since $X=\overline{M}$ (no point of $X$ is isolated), $X=\bigcup_\alpha P_\alpha$. Conversely, if $X$ is a finite union of affine planes, then $K\equiv0$ and equality holds by Theorem~\ref{thm:main}.
\end{proof}

\begin{proof}[Proof of Corollary~\ref{cor:complex}]
A complex algebraic curve $X\subset\C^n=\R^{2n}$ is a closed semialgebraic set. Its complex singular set $\Sing^{\C}(X)$ is finite, and $X\setminus\Sing^{\C}(X)$ is a real analytic surface, so $X$ is pure two-dimensional and $\Sing_2(X)\subset\Sing^{\C}(X)$. Conversely, if $p\in X$ is a point at which $X$ is a $C^1$ submanifold, then $\theta(X,p)=1$ (Subsection~\ref{subsec:density}), hence $m(X,p)=1$ by Thie's theorem \cite{Thie:1967}, so $p$ is a smooth point of $X$. Therefore $\Sing_2(X)=\Sing_1(X)=\Sing^{\C}(X)=\{p_1,\ldots,p_s\}$. By \cite{Thie:1967} we have $\theta(X,p_i)=m(X,p_i)$, and $\theta(X,\infty)=\deg(X)$ by \cite{Stoll:1964} (see also \cite{SilvaSampaio:2026}). The formula follows from Theorem~\ref{thm:main}.

The regular part of a complex curve is a minimal surface (the second fundamental form satisfies $\II(Jv,Jw)=-\II(v,w)$, where $J$ is the complex structure, so that $\II(e_1,e_1)+\II(Je_1,Je_1)=0$). Hence Corollary~\ref{cor:minimal} gives the inequality, with equality if and only if $X$ is a finite union of affine real $2$-planes. Each of these planes is the closure of a connected component of $\Reg(X)$, hence it is a complex curve, i.e., a complex affine line; and a union of $d$ distinct complex lines has degree $d$.
\end{proof}

\bigskip

\noindent {\bf Acknowledgments:} 
During the preparation of this work, the authors used the AI assistants Claude, Gemini, and ChatGPT to refine arguments, locate references, and draft and revise the manuscript. All AI-generated suggestions were independently verified by the authors, who take full responsibility for the content and accuracy of the final manuscript.

\end{document}